\documentclass[12pt]{amsart}
\usepackage{amscd}
\usepackage{latexsym,amssymb}

\newtheorem{theorem}{Theorem}[section]
\newtheorem{corol}[theorem]{Corollary}
\newtheorem{prop}[theorem]{Proposition}
\newtheorem{conj}[theorem]{Conjecture}
\theoremstyle{definition}
\newtheorem{defin}[theorem]{Definition}

\theoremstyle{remark}

\newtheorem*{erem}{Remark}

\def\DMO{\DeclareMathOperator}
\def\md#1#2#3{#1\equiv#2\pmod{#3}}
\def\cf#1{(cf.~\cite{#1})}  \DMO{\Aut}{Aut}
\def\mZ{\mathbb Z}	\def\mA{\mathbb A}
\def\sA{\mathsf A}	\def\sB{\mathsf B}
\DMO{\im}{Im}		\def\iso{\simeq}
\def\+{\oplus}		\def\bop{\bigoplus}
\def\al{\alpha} \def\be{\beta} \def\ga{\gamma}
\DMO{\add}{add}	\DMO{\Mat}{Mat}
\DMO{\spec}{Spec}	
	
\def\mod#1#2{{#2}\mbox{-}\mathsf{mod}^{#1}}
\def\fab{\mathsf{fab}}	
\def\oP{\overline P}	\def\op{\overline p}
\def\oh{\overline h}	
\DMO{\prd}{pr.dim}	\def\gai#1{\stackrel{\ga_#1}\larr}
\def\kk#1{{}\!#1\!{}}	\def\ou{\ol u}
\def\bR{\mathbf R}	\def\bA{\mathbf A}	\def\bB{\mathbf B}
\def\bD{\mathbf D}	\def\bC{\mathbf C}
\def\kM{\mathcal M}  \def\kN{\mathcal N}  \def\kC{\mathcal C}
\def\vi{\varphi}   \def\Si{\Sigma}   \def\xx{\times}
\def\iff{if and only if }  \def\mps{\mapsto}  \def\ol{\overline}
\def\Id{\mathrm{Id}}	\def\oM{\overline M}   \DMO{\End}{End}
\DMO{\Hom}{Hom}	\DMO{\id}{id}
\def\gnr#1{\left\langle\,#1\,\right\rangle}
\def\set#1{\left\{\,#1\,\right\}}  \def\={\setminus}
\def\rap{\noindent} 
\def\lst#1#2{#1_1,#1_2,\dots,#1_{#2}}
\def\mat#1{\begin{matrix}#1\end{matrix}}
\def\mtr#1{\begin{pmatrix}#1\end{pmatrix}}
\def\oC{\overline\bC}	\def\Mod#1{{#1}\mbox{-}\mathsf{Mod}} 
\def\oc{one-to-one correspondance}
\def\setsuch#1#2{\left\{\,#1\,|\,#2\,\right\}}
\def\mQ{\mathbb Q} \def\mN{\mathbb N}  \def\nT{\mathrm t}
\def\La{\Lambda}  \def\la{\lambda}  \def\si{\sigma}  \def\Ga{\Gamma}
\def\th{\theta}  \def\om{\omega} \DMO{\chr}{char}
\def\bE{\mathbf E} \def\fG{\mathbf g}  \def\fK{\mathbf k}
\def\larr{\longleftarrow} \def\8{\infty}  \def\sbe{\subseteq}
\def\rB{\mathrm B}  \def\rD{\mathrm D}  \def\rS{\mathrm S}
\def\sC{\mathsf C}  \def\sQ{\mathsf Q}  \def\sL{\mathsf L}
\def\dE{\mathfrak E}  \def\dF{\mathfrak F}  \def\dX{\mathfrak X}
\def\de{\delta}  \def\bH{\mathbf H}

\title{On cubic functors}
\author{Yuriy A. Drozd}

\thanks{This work was supported by the Max-Planck-Institut f\"ur Mathematik 
and by the CRDF Award UM2-2094. }

\begin{document}
\maketitle

\begin{center}
 Kiev Taras Shevchenko University\\ Department of Mechanics and Mathematics\\
01033 Kiev, Ukraine\\ E-mail: {\it yuriy@drozd.org}
\end{center}

\vskip.5in
{\bf Abstract}.
 We prove that the description of cubic functors is a wild problem in the sense
 of the representation theory. On the contrary, we describe several special classes
 of such functors (\emph{2}-divisible, weakly alternative, vector spaces and torsion
 free ones). We also prove that cubic functors can be defined locally and obtain
 corollaries about their projective dimensions and torsion free parts.

\tableofcontents

\bigskip
\section*{Introduction}

 Polynomial functors appeared in algebraic topology (\cf{em}) and
 proved to be useful for lots of questions in homotopy theory. That is
 why their study seems to be of interest. Mostly, one deals with
 ``continuous'' polynomial functors, which are defined by their values
 on free groups, and in what follows we consider only such ones. In
 \cite{qu} the author gave a description of all (finitely generated)
 quadratic functors. This description was done using the technique of
 the so called \emph{matrix problems}, namely, representations \emph
 { of bunches of chains}. The reduction to such a matrix problem
 depends on the presentation of quadratic functors as modules over a
 special ring, which happens to be an order in a semi-simple algebra.
 Then we used the procedure of \cite{pure}. Now we are going to pass
 the same way for \emph{cubic functors}. We check that the
 corresponding ring is again an order in a semi-simple algebra.
 Unfortunately, in this case, the classification of modules is a
 \emph{wild problem}, i.e., it includes, in some sense, the
 classification of all representations of all finitely generated
 algebras over the residue field $\,\mZ/2\,$
 (cf. Section~\ref{sec2}). Thus, there can be no ``good'' description
 of all cubic functors.

 On the contrary, such a description becomes possible (and rather
 analogous to that of quadratic functors) if we ``make $2$
 invertible,'' i.e., consider cubic modules over the ring
 $\,\mZ[\,1/2\,]\,$. Such a description is given in
 Section~\ref{sec3}. Just as for quadratic case, this classification
 problem is \emph{tame}, i.e., indecomposable modules depend on some
 ``discrete'' combinatorial data and on at most one ``continuous
 parametre,'' which is an irreducible polynomial from
 $\,\mZ/2[\,t\,]\,$ (or, the same, a closed point of
 $\,\mA^1_{\mZ/2}=\spec\mZ/2[\,t\,]\,$).  We also consider
 \emph{weakly alternative} cubic functors $\,F\,$, i.e., those with
 $\,F(\mZ)=0\,$. Their classification given in Section~\ref{sec4} is
 also tame, though this time the corresponding ring has both torsion
 and nilpotent  ideals.  At last, we give a classification of cubic functors with the
 image being \emph{vector spaces} (Section~\ref{sec5}) as well as
 of  \emph{torsion free} ones (Section~\ref{sec6}). These problems also
 happens to  be tame.  We end up with one conjecture concerning
 polynomial functors of higher degrees that arises from the parallel between
 quadratic and \emph{2}-divisible cubic functors and some corollaries from this
 conjecture. 
 
\section{Generalities}
\label{sec1}

 We suppose all categories to be \emph{pre-additive}, i.e., such that
 their morphism sets are abelian groups and the composition is
 bi-additive. On the contrary, we do not suppose the \emph{functors} to be
 additive, though we always suppose that $\,F(0)=0\,$ for a zero
 morphism. Remind that a \emph{fully additive} category is an additive
 category such that every idempotent in it splits. If
 $\,F:\sA\to\sB\,$ is a functor from an additive category $\,\sA\,$ to
 a fully additive category $\,\sB\,$ and
 $\,A=\bop_{k=1}^nA_k\,$ is an object from $\,\sA\,$, consider the
 corresponding embeddings $\,i_k:A_k\to A\,$ and projections
 $\,p_k:A\to A_k\,$. Put $\,e_k=i_kp_k\,$ and $\,f(k)=F(e_k)\,$.
 Certainly, $\,e_k\,$, hence, $\,f(k)\,$, are pairwise orthogonal
 idempotents and $\,\im f(k)\iso F(A_k)\,$. Moreover, put
 $\,f(kl)=F(e_k+e_l)-f(k)-f(l)\,$ ($\,k<l\,$). Then $\,f(kl)\,$ are
 also idempotents, pairwise orthogonal and orthogonal to all
 $\,f(k)\,$. Hence, $\,F(A)\,$ has a direct summand
 $\,(\bop_kF(A_k))\+(\bop_{k<l}F_2(A_k|A_l))\,$, where one denotes
 $\,F_2(A_k|A_l)=\im f(kl)\,$. Define recursively for any $\,1\le
 k_1<k_2<\dots<k_r\le n\,$
$$
  f(k_1k_2\dots k_m)=F(e_{k_1}+\dots e_{k_m})-
  \sum_{r<m}\sum_{1\le l_1<\dots<l_r\le m}f(k_{l_1}k_{l_2}\dots
  k_{l_r}) 
$$
 and $\,F_m(A_1|A_2|\dots|A_m)=\im f(k_1k_2\dots k_m)\,$. Then
$$
 F(A)=\bop_{m\le n}\left(\bop_{k_1k_2\dots k_m}
 F_m(A_{k_1}|A_{k_2}|\dots|A_{k_m}) \right)\,. 
$$
 One easily sees that $\,F_m\,$ is indeed a functor
 $\,\sA^m\to\sB\,$. It is called the $\,m$-th \emph{cross-effect} of
 $\,F\,$. Of course, $\,F_1=F\,$. If there is a positive integer
 $\,n\,$ such that $\,F_m=0\,$ for $\,m>n\,$, one calls $\,F\,$ a
 \emph{polynomial functor}. The smallest possible $\,n\,$ with this
 property is called the \emph{degree} of $\,F\,$. Certainly, the
 functors of degree $1$ are just additive ones. The functors of degree
 $1$ are called \emph{linear}, of degree $2$ \emph{quadratic}, of
 degree $3$ \emph{cubic}, etc. For an arbitrary pre-additive category
 $\,\sA\,$, one can always consider its \emph{fully additive hull}
 $\,\add\sA\,$ and we call polynomial functors $\,\sA\to\sB\,$ those
 from $\,\add\sA\,$ to $\,\sB\,$. Of course, for additive functors we
 get in this way nothing new as every additive functor $\,\sA\to\sB\,$
 can be uniquely (up to isomorphism) prolonged to an additive functor
 $\,\add\sA\to\sB\,$.

 In what follows, we consider the case when $\,\sA=\mZ\,$, the ring of
 integers, and $\,\sB=\Mod\bR\,$, the category of modules over a ring
 $\,\bR\,$, mainly even $\,\sB=\mod{}\bR\,$, the category of finitely
 generated $\,\bR$-modules. Note that $\,\add\mZ=\fab\,$, the
 category of finitely generated free abelian groups. In this case
 polynomial functors are 
 called \emph{polynomial $\,\bR$-modules} and the category of finitely
 generated polynomial $\,\bR$-modules of degree at most $\,n\,$ is
 denoted by $\,\mod n\bR\,$. Of course, $\,\mod1\bR\,$ coincide with
 the category of finitely generated $\,\bR$-modules. If $\,\bR=\mZ\,$,
 we are just speaking of ``polynomial modules'' not precising the target
 category. 

 To define a polynomial functor $\,F:\fab\to\sB\,$, one only has to
 define the following data:
\begin{itemize}
 \item
  objects $\,F_m=F_m(\mZ|\mZ|\dots|\mZ)\,$;
 \item
  morphisms $\,h^m_k:F_m\to F_{m+1}\,$ that are compositions of the
  morphisms 
\begin{align*}
 F_m(A_1|A_2|\dots|A_m)&\to F(A_1\+A_2\+\dots\+A_m)\\&\to
 F(A_1\+\dots\+A_k\+A_k\+\dots\+A_m)\\&\to
 F_{m+1}(A_1|\dots|A_k|A_k|\dots|A_m)\,, 
\end{align*}
 where all $\,A_j=\mZ\,$, the first morphism is the embedding of the
 direct summand, the last one is the projection onto the direct
 summand and the middle one is $\,F(\ga)\,$, $\,\ga|_{A_j}\,$ being
 identity if $\,j\ne k\,$ and $\,\ga|_{A_k}\,$ being the diagonal 
 embedding $\,A_k\to A_k\+A_k\,$;
 \item
  morphisms $\,p^m_k:F_{m+1}\to F_m\,$ that are compositions of the
 morphisms
\begin{align*}
 F_{m+1}(A_1|\dots|A_k|A_k|\dots|A_m)&\to
 F(A_1\+\dots\+A_k\+A_k\+\dots\+A_m)\\&\to F(A_1\+A_2\+\dots\+A_m)
 \\&\to F_m(A_1|A_2|\dots|A_m)\,,
\end{align*}
 where again all $\,A_j=\mZ\,$, the first morphism is the embedding of
 the direct 
 summand, the last one is the projection onto the direct summand and
 the middle one is $\,F(\be)\,$, $\,\be|_{A_j}\,$ being identity for
 $\,j\ne k\,$, while $\,\be|_{A_k}\,$ being the codiagonal (summation)
 mapping $\,A_k\+A_k\to A_k\,$.
\end{itemize}

 In particular, a cubic $\,\bR$-module $\,F\,$, where $\,\bR\,$ is a
 ring, is defined by a diagram of $\,\bR$-modules:
\begin{equation}\label{eq2}
\begin{CD}
F_1 \mat{@>h>>\\@<p<<} F_2
 \mat{@>{h_1}>>\\@>{h_2}>>\\@<{p_1}<<\\@<{p_2}<<}
 F_3 
\end{CD}
\end{equation}
 One can show that such a diagram corresponds to a cubic module \iff
 the following relations hold:
\begin{equation}\label{eq1}
\begin{split}
 &  h_1p_2=h_2p_1=0\,,\quad h_1h=h_2h\,,\quad pp_1=pp_2\,,\\
 &  h_ip_ih_i=2h_i\,,\quad p_ih_ip_i=2p_i\quad (i=1,2)\,,\\
 &  hph=2(h+(p_1+p_2)\ol h)\,,\quad php=2(p+\ol p(h_1+h_2))\,,\\
 &  \ol hp+h_1+h_2=h_1p_1h_2p_2h_1+h_2p_2h_1p_1h_2\,,\\
 &  h\ol p+p_1+p_2=p_1h_2p_2h_1p_1+p_2h_1p_1h_2p_2\,,
\end{split}
\end{equation}
 where $\,\ol h=h_1h=h_2h\,$ and $\,\ol p=pp_1=pp_2\,$ \cite{dre}. 

 Consider the pre-additive category $\,\bA\,$ with three objects
 $\,1,2,3\,$ and generating morphisms
$$
 h:1\to2\,,\ p:2\to1\,,\ h_i:2\to3\,,\ p_i:3\to2\quad (i=1,2)
$$
 subject to the relations \eqref{eq1}. Then cubic $\,\bR$-modules are
 just linear functors $\,\bA\to\mod{}\bR\,$. Sometimes we identify
 $\,\bA\,$ with the endomorphism ring $\,\End(1\+2\+3)\,$. Obviously,
 the categories of modules over $\,\bA\,$ and over this ring coincide.

\section{Wildness}
\label{sec2}

 We are going to show that the description of \emph{all} cubic
 functors is a \emph{wild problem} in the sense of the representation
 theory \cf{tw}. Moreover, we show that it is even true for those
 cubic functors ``freely generated by their 1-part,'' or, the same,
 for $\,\bA(1,1)$-modules.

\begin{prop}
 The ring $\,\bA(1,1)\,$ is isomorphic to the subring of $\,\mZ^3\,$
 with the basis $\,\set{(1,1,1),(2,0,0),(0,6,0)}\,$.
\end{prop}
\begin{proof}
 Put $\,a=ph-\ol{ph},\,b=ph\,$. One easily verifies that these two
 elements generate $\,\bA(1,1)\,$, $\,ab=ba=0\,$, $\,a^2=2a\,$ and
 $\,b^2=6b\,$. Hence, one gets an isomorphism  by mapping
 $\,a\mps(2,0,0)\,$, $\,b\mps(0,6,0)\,$.
\end{proof}

 The category of $\,\bA(1,1)$-modules can be embedded into that of
 $\,\bA$-modules. Namely, every $\,\bA(1,1)$-module $\,M\,$ gives rise
 to the cubic module $\,F\,$, where
$$
 F_1=M\,,\quad F_2=\bA(1,2)\*_{\bA(1,1)}M\quad\textrm{and}\quad
 F_3=\bA(1,3)\*_{\bA(1,1)}M\,
$$
 with the obvious action of morphisms.
 It is known (and easy to check) that this procedure really defines a
 functor, which is a full embedding.

 Let $\,\Si_n\,$ be the free (non-commutative) algebra
 $\,\mZ/4\gnr{\lst xn}\,$ with $\,n\,$ generators over the residue
 ring $\,\mZ/4\,$.

\begin{prop}\label{prop2}
 Denote by $\,\vi:\bA(1,1)\to\Si_2\,$ the homomorphism mapping $\,a\,$
 to $\,2x_1\,$ and $\,b\,$ to $\,2x_2\,$. For any $\,\Si_2$-module
 $\,L\,$, let $\,^\vi L\,$ be the $\,\bA(1,1)$-module obtained from
 $\,L\,$ via the ``change of rings'' $\,\vi\,$. Then, for any
 $\,\Si_2$-modules $\,L,L'\,$ which are free as $\,\mZ/4$-modules,
\begin{itemize}
 \item
 $\,{}^\vi L\iso{}^\vi L'\,$ as $\,\bA(1,1)$-modules \iff $\,L/2L\iso
 L'/2L'\,$ as $\,\Si_2$-modules; 
 \item
 $\,\bA(1,1)$-module $\,{}^\vi L\,$ is indecomposable \iff so is
 $\,\Si_2$-module $\,L/2L\,$.
\end{itemize}
\end{prop}
\begin{proof}
 [\rm The proof is evident.]
\end{proof}

\begin{corol}
 For every $\,n\,$ there is a cubic module $\,\kM\in\kC(\Si_n)\,$
 such that, for any $\,\Si_n$-modules $\,L,L'\,$ which are free over
 $\,\mZ/4\,$,
\begin{itemize}
 \item
  $\,\kM\*_{\Si_n}L\iso\kM\*_{\Si_n}L'\,$ \iff $\,L/2\iso L'/2\,$;
 \item
  $\,\kM\*_{\Si_n}L\,$ is indecomposable \iff $\,L/2\,$ is
  indecomposable. 
\end{itemize}
\end{corol}
\begin{proof}
 One only has to consider the $\,\Si_2$-$\,\Si_n$-bimodule $\,\kN\,$
 which is free of rank $\,n+1\,$ over $\,\Si_n\,$, while the action of
 $\,\Si_2\,$ is given by the following rules:
\begin{itemize}
 \item
  $\,x_1\,$ is acting as the (upper) Jordan cell;
 \item
  $\,x_2\,$ is acting as the matrix
$$
\mtr{0&0&\dots&0&0\\
     x_1&0&\dots&0&0\\
     0&x_2&\dots&0&0\\
     \hdotsfor{5}\\
     0&0&\dots&x_n&0}
$$
\end{itemize}
 Afterwards, one puts $\,\kM={}^\vi\kN\,$ with $\,\vi\,$ defined
 as in Proposition \ref{prop2}. All properties are
 verified immediately (cf., e.g., \cite{tw}).
\end{proof}

 This Corollary shows that the classification of all cubic modules
 (even of $\,\bA(1,1)$-modules) is at least so
 complicated as the classification of $\,n$-tuples of matrices over
 the filed $\,\mZ/2\,$ up to conjugations, for all $\,n\,$. It is
 just the situation when one calls such a problem ``wild.'' 

\section{\emph{2}-divisible case}
\label{sec3}

 In this section, we consider the case of cubic \emph{2-divisible}
 modules, i.e., cubic functors $\,F:\fab\to\mod{}{\mZ'}\,$, where
 $\,\mZ'=\mZ[1/2]\,$. It happens that here we are able to describe all
 (finitely generated) cubic modules. Moreover, this description is
 quite alike that of quadratic $\,\mZ$-modules given in
 \cite{qu}. Certainly, \emph2-divisible cubic modules correspond to the
 modules over the category $\,\bA'=\bA[1/2]\,$.

\begin{prop}\label{prop31}
 The category $\,\bA'\,$ is Morita equivalent to the direct product
 $\,\mZ'\xx\mZ'\xx\bB\,$ where $\,\bB\,$ is the subring in
 $\,\mZ'\xx\Mat(2,\mZ')\xx\Mat(2,\mZ')\xx\mZ'\,$ consisting of all
 quadruples of the form:
\begin{align*}
&\left(a,\mtr{b_1&3b_2\\b_3&b_4},\mtr{c_1&3c_2\\c_3&c_4},d\right)\\
\intertext{with}
&a\equiv b_1,\ b_4\equiv c_1,\ c_4\equiv d \pmod3
\end{align*}
\end{prop}
\begin{erem}
 One can see that a module over the ring $\,\bB\,$ is given by a
 diagram of abelian (\emph2-divisible) groups of the following shape: 
\begin{equation}\label{eq1d}
\begin{CD}
M_1 \mat{@>\al_1>>\\@<\be_1<<} M_2 \mat{@>\al_2>>\\@<\be_2<<}
 M_3 
\end{CD},
\end{equation}
 satisfying the relations:
\begin{equation*}
\begin{split}
&  \al_1\be_1\al_1=3\al_1\,,\quad \be_1\al_1\be_1=3\be_1\,,\\
&  \al_2\be_2\al_2=3\al_2\,,\quad \be_2\al_2\be_2=3\be_2\,,\\
&   \be_1\be_2=0\,,\quad  \al_2\al_1=0\,,\\
&  \al_1\be_1+\be_2\al_2=3\id_2\,.
\end{split}
\end{equation*}
 We always identify $\,\bB$-modules with such diagrams.
\end{erem}
\begin{proof}
 Put $\,\bA_i=\bA'(i,i)\,$ ($\,i=1,2,3\,$) and consider the
 following elements:
\begin{align*}
\intertext{in $\,\bA_1\,$:}
& f_1=\frac{ph-\op\oh}2,\ e_1=\id_1-f_1,\ a_1=\frac{\op\oh}2;\\
\intertext{in $\,\bA_2\,$:}
& e_2=\frac{p_1h_1}2,\ f_2=\frac{p_2h_2}2,\ g=\id_2-e_2-f_2,\\
& u=hp,\ g_1=\frac{gu}2,\ g_2=g-g_1,\\
& v_i=p_ih_j\ (i,j=1,2,\,i\ne j),\ v_1'=\frac{5v_1-v_1v_2v_1}2,\\
& a_2=v_1v_2,\ b_2=v_2v_1;\\
\intertext{in $\,\bA_3\,$:}
& u_i=h_ip_i\ (i=1,2),\ e_3=\frac{u_1}2,\ f_3=\id_3-e_3,\\
& \ou=\oh\op+e_3,\ a_3=2u_2e_3-\ou,\ b_3=2e_3u_2-\ou;\\
\intertext{in $\,\bA'(2,1)\,$:}
& p'=p-\frac{\op(h_1+h_2)}2.
\end{align*}
 On can verify that:
\begin{itemize}
 \item
  elements $e_i,\,f_i,\,g_i\,$ are orthogonal idempotents;
 \item
  elements $\,e_1,f_1,a_1\,$ form a $\,\mZ'$-basis of $\,\bA_1\,$;
 \item
  $a_1=e_1a_1e_1\,$ and $\, a_1^2=3a_1$\,;
 \item
  elements $\,e_2,f_2,g_1,g_2,a_2,b_2,v_1,v_2,v_1b_2,a_2v_1\,$ form a
  $\,\mZ'$-basis of $\,\bA_2\,$;
 \item
  $\,g\bA_2=\bA_2g=\gnr{g_1,g_2}\,$;
 \item
  $\,a_2=e_2a_2e_2\,$, $\,b_2=f_2b_2f_2\,$, $\,a_2^2=3a_2\,$,
  $\,b_2^2=3b_2\,$; 
 \item
  $\,v_1'(v_2/2)=e_2\,$, $\,(v_2/2)v_1'=f_2\,$; 
 \item
  $\,p'(h/2)=f_1\,$, $\, (h/2)p'=g_1\,$; 
 \item
  elements $\,e_3,f_3,a_3,b_3,a_3b_3,b_3a_3\,$ form a $\,\mZ'$-basis of
  $\,\bA_3\,$; 
 \item
  $\,a_3=f_3a_3e_3\,$, $\,b_3=e_3b_3f_3\,$, $\,a_3b_3a_3=3a_3\,$,
  $\,b_3a_3b_3=3b_3\,$; 
 \item
  $\,p'(h/2)=f_1\,$ and $\,(h/2)p'=g_1\,$;
 \item
  $\,(p_1/2)h_1=e_2\,$ and  $\,h_1(p_1/2)=e_3\,$;  
 \item
  $\,gp_i=h_ig=0\,$ ($\,i=1,2\,$) and $\,g_2h=pg_2=0\,$;
 \item
  $\,f_2h_i=p_if_2=0\,$ ($\,i=1,2\,$).
\end{itemize}
 Consider the projective $\,\bA'$-modules:
\begin{align*}
& E_1=\bA'(1,\_)e_1,\ F_1=\bA'(1,\_)f_1,\\
& E_2=\bA'(2,\_)e_2,\ F_2=\bA'(2,\_)f_2,\ G_i=\bA'(2,\_)g_i,\\ 
& E_3=\bA'(3,\_)e_3,\ F_3=\bA'(3,\_)f_3,\\
\end{align*}
 They are projective generators of the category of
 $\,\bA'$-modules. The above equalities imply that:
\begin{align*}
& F_1\iso G_1,\ E_2\iso F_2,\ E_2\iso E_3;\\
& \Hom(G_i,G_j)=0\ \textrm{ if }\ i\ne j;\\
& \Hom(G_i,E_j)=\Hom(E_j,G_i)=\Hom(G_i,F_3)=\Hom(F_3,G_i)=0;\\
& \Hom(E_1,F_3)=\Hom(F_3,E_1)=0;\\
& \End E_1\iso e_1\bA_1e_1=\gnr{e_1,a_1}\iso\bD;\\
& \End E_2\iso e_2\bA_2e_2=\gnr{e_2,a_2}\iso\bD;\\
& \End F_3\iso f_3\bA_3f_3=\gnr{f_3,a_3b_3}\iso\bD;\\
& \End G_i=\mZ';\\
& \Hom(E_3,F_3)=\gnr{b_3}\iso\mZ';\\
%\end{align*}
%\begin{align*}
& \Hom(F_3,E_3)=\gnr{a_3}\iso\mZ';\\
& \Hom(E_1,E_3)=\gnr{\op}\iso\mZ';\\
& \Hom(E_3,E_1)=\gnr{\oh}\iso\mZ'.\hfill
\end{align*}
 Here $\,\bD\,$ denotes the ring $\,\mZ'[t]/(t^2-3)\,$, which is
 isomorphic to the subring in $\,\mZ'\xx\mZ'\,$ consisting of the
 pairs $\,\setsuch{(a,b)}{\md a b3}\,$. Therefore, the category
 of $\,\bA'$-modules is equivalent to that of $\,\bE$-modules, where
 $\,\bE=\End (G_1\+G_2\+E_1\+F_3\+E_3)\iso \mZ'\xx\mZ'\xx\bB\,$.
\end{proof}

 One can verify that the $\,\bA'$-modules $\,G_i\,$ correspond to
 the functors $\,\rS^2\,$ and $\,\La^2\,$, where $\,\rS^r\,$
 and $\,\La^r\,$ denote, as usually, the $\,r$-th symmetric and
 exterior powers. Actually, these functor are quadratic, and it is the
 main reason that they ``stand apart'' really cubic ones. 

\medskip
 Now the theory of \emph2-divisible cubic modules becomes quite
 parallel to that of quadratic modules \cite{qu} and one gets the main
 results just following the same way. As the proofs are also
 almost the same, we replace them by the exact
 references to the corresponding items from \cite{qu}. We denote by
 $\,e_i\,$ ($\,i=1,2,3\,$) the natural idempotents in $\,\bB\,$:
\begin{align*}
 e_1&= \left(1,\mtr{1&0\\0&0},\mtr{0&0\\0&0},0 \right)\,,\\
 e_2&= \left(0,\mtr{0&0\\0&1},\mtr{1&0\\0&0},0 \right)\,,\\
 e_1&= \left(0,\mtr{0&0\\0&0},\mtr{0&0\\0&1},1 \right)\,.
\end{align*}
 Put $\,\sB_i=\bB e_i\,$; they are just indecomposable projective
 $\,\bB$-modules. Denote by $\,\sL_1\,$ the projection of $\,\sB_1\,$
 onto the first component, by $\,\sL_2\,$ its projection onto the
 second component; by $\,\sL_3\,$ the projection of $\,\sB_2\,$ onto
 the second component, by $\,\sL_4\,$ its projection onto the third
 component; by $\,\sL_5\,$ the projection of $\,\sB_3\,$ onto the
 third and by $\,\sL_6\,$ its projection onto the forth
 component. Then $\,\sL_i\,$ are just irreducible torsion free
 $\,\sB$-modules. Put $\,\mZ_{(p)}=\setsuch{a/b}{a,b\in\mZ,\,p\nmid
 b}\,$, $\,M_{(p)}=\mZ_{(p)}\*M\,$.

\begin{theorem}\label{t21}
\begin{enumerate}
 \item 
 2-divisible cubic modules $\,M,N\,$ are isomorphic \iff
 $\,M_{(p)}\iso N_{(p)}\,$ for every odd prime $\,p\,$.   
\item 
 Given a set $\,\set{N_p}\,$ ($\,p>2\,$) where $\,N_p\,$ is a
 cubic $\,\mZ_{(p)}$-module, there is a 2-divisible cubic module
  $\,M\,$ such that $\,M_{(p)}\iso N_p\,$ for all $\,p\,$ \iff almost
  all (i.e., all but a finite set) of them are torsion free (maybe,
  zero) and $\,\mQ\*N_p\iso\mQ\*N_q\,$ for all $\,p,q\,$. In this case,
  $\,\mQ\*M=\mQ\*N_p\,$ and $\,\nT M=\bop_p\nT N_p\,$.
\end{enumerate}
\end{theorem}
\begin{proof}
 [\rm Cf] \cite[Theorem 2.1]{qu}.
\end{proof}

 Now we define our main personages: string and band modules. We
 rearrange a bit the description given in \cite[Section 4]{qu} to make
 it more comprehensible. We start with some notations and definitions.

\begin{defin}
\begin{enumerate}
 \item
  Define an equivalence  relation $\,-\,$ on the set $\,\set{1,2,3,4,5,6}\,$ 
 with the only non-trivial equivalences $\,2-3\,$
  and $\,4-5\,$ and a symmetric relation $\,\sim\,$ (not equivalence!)
  putting $\,1\sim2,\,3\sim4,\,5\sim6\,$. 
 \item
  Define the following mappings acting on every diagram of the
  form \eqref{eq1d}:
\begin{equation}\label{theta}
\begin{array}{lp{1cm}l}
 \th(11)=3\id_1-\be_1\al_1\,,&&\th(22)=\be_1\al_1\,,\\
 \th(23)=\al_1\,,&&\th(32)=\be_1\,,\\
 \th(33)=\be_2\al_2\,,&&\th(44)=\al_2\be_2\,\\
 \th(45)=\al_2\,,&&\th(54)=\be_2\,,\\
 \th(55)=\al_2\be_2\,&&\th(66)=3\id_3-\al_2\be_2\,.
\end{array}
\end{equation}
  \item
  Define the function
  $\,\nu:\set{\set{1,2},\set{3,4},\set{5,6}}\to\mN\,$ putting
  $\,\nu\set{1,2}=1,\ \nu\set{3,4}=2,\ \nu\set{5,6}=3 \,$.
  \item
  A \emph{primary polynomial} over a field $\,\fK\,$ is, by
  definition, a power of an irreducible polynomial with the leading
  coefficient $\,1\,$.
\end{enumerate}
\end{defin}

 Now the string and band modules can be defined as follows.

\begin{defin}\label{str1}
 \begin{enumerate}
  \item
  A \emph{string diagram} $\,\rD\,$ is a diagram of one of the
  following types: 
\begin{align*}
 \begin{array}{*{11}c}
 j_1& & & &j_2j_3& & & &j_{2n-2}j_{2n-1}& & \\
 &\kk{k_1}& &\kk{k_2}& &\kk{k_3}& &\kk{\cdots}&
 &\kk{k_{2n-1}} \\ 
 & & i_1i_2& & & &i_3i_4& & & &i_{2n-1}
 \end{array},
 \tag{i}\\\intertext{or}
 \begin{array}{*{11}c}
  & &j_2j_3& & & &j_{2n-2}j_{2n-1}& & \\
  &\kk{k_2}& &\kk{k_3}& &\kk{\cdots}&
 &\kk{k_{2n-1}} \\ 
  i_2& & & &i_3i_4& & & &i_{2n-1}
\end{array}, \hspace*{2em}
 \tag{ii}\\\intertext{or}
 \begin{array}{*{11}c}
 j_1& & & &j_2j_3&  & & &j_{2n} \\
 &\kk{k_1}& &\kk{k_2}& &\kk{\cdots}&
 &\kk{k_{2n}}& \\ 
 & & i_1i_2& & & &i_{2n-1}i_{2n}& &
 \end{array}, \hspace*{3em}
 \tag{iii}
\end{align*}
 where $\,k\in\mN\,$, $\,i_m,j_m\in\set{1,2,3,4,5,6}\,$, 
 satisfying the following conditions:
 \begin{itemize}
  \item
  $\,i_{2m-1}\sim i_{2m}\,$ for every $\,m=1,\dots,n\,$.\\
 (This condition is empty for diagrams (i) and (ii) if $\,m=n\,$ and
 for diagrams (ii) if $\,m=1\,$; nevertheless, in these cases we
 \emph{define} $\,i_1\,$ or $\,i_{2n}\,$ so that this condition
 holds.) 
  \item
  $\,j_{2m+1}\sim j_{2m}\,$ for every $\,m=1,\dots,n-1\,$.
  \item
  $\,i_k-j_k\,$ for every $\,m=1,\dots,2n\,$.\\
 (This condition is empty for diagrams (i) and (ii) if $\,m=2n\,$ and
 for diagrams (ii) if $\,m=1\,$.)
 \end{itemize}
  \item
  The \emph{string $\,\bB$-module} corresponding to a string diagram
  $\,\rD\,$ is the $\,\bB$-module $\,M=M^\rD\,$ with the generators 
$$
 \lst\fG n\,,\qquad \fG_m\in M_{\nu\set{i_{2m-1},i_{2m}}}\,,
$$
 and the defining relations
$$
 3^{k_{2m}}\th(i_{2m}j_{2m})g_m =
 3^{k_{2m+1}}\th(i_{2m+1}j_{2m+1})g_{m+1}
$$
 for $\,i=0,\dots,n\,$; we put here $\,g_0=g_{n+1}=0\,$ and omit the
  relation with $\,m=n\,$ for diagrams (i), (ii) and with
  $\,m=1\,$ for diagrams (ii).
  \item
  A \emph{string cubic module} is one corresponding to a string
 $\,\bB$-module via the Morita equivalence of
 Proposition~\ref{prop31}. 
 \end{enumerate}
\end{defin}

\begin{defin}
\begin{enumerate}
 \item
 A  diagram $\,\rD\,$ of type (iii) is said to be
 \emph{non-periodic} if it cannot be obtained by a repetition of a
 smaller diagram $\,\rD_1\,$ of the same type:
 $\,\rD\ne\rD_1\rD_1\dots\rD_1\,$. 
 \item
 A \emph{band data} is a pair $\,\rB=(\rD,f(t))\,$ consisting of a
 non-periodic diagram of type (iii) with the additional
 condition:
 \begin{itemize}
  \item
  $\,j_{2n}\sim j_1\,$,
 \end{itemize}
 and of a primary polynomial
 $\,f(t)=\la_1+\la_2t+\dots+\la_dt^{d-1}+t^d\,$ over the residue field
 $\,\mZ/3\,$ with $\,\la_1\ne0\,$.
 \item
  The \emph{band $\,\bB$-module} corresponding to a band data
 $\,\rB\,$ is the $\,\bB$-module $\,M=M^\rB\,$ with the generators
$$
\hfil \fG_{ml} \ \ (m=1,\dots,n,\,l=1,\dots,d)\,,
 \ \ \fG_{ml}\in M_{\nu\set{i_{2m-1},i_{2m}}}\,,
$$
 and the defining relations:
\begin{align*}
 3^{k_{2m}}\th(i_{2m}j_{2m})g_{ml} &=
 3^{k_{2m+1}}\th(i_{2m+1}j_{2m+1})g_{m+1,l}\quad \mathrm{ if }\quad 
 1\le m<n \,;\\ 
 3^{k_{2n}}\th(i_{2n}j_{2n})g_{nl} &=
 3^{k_1}\th(i_1j_1)g_{1,l+1}\quad \mathrm{ if }\quad 1\le l<d \,;\\ 
 3^{k_{2n}}\th(i_{2n}j_{2n})g_{nd} &=
 -3^{k_1}\th(i_1j_1)\sum_{\nu=1}^d\la_\nu g_{1\nu}\,.
\end{align*}
 \item
  A \emph{band cubic module} is one corresponding to a band
  $\,\bB$-module via the Morita equivalence of
  Proposition~\ref{prop31}. 
\end{enumerate}
\end{defin}
\rap
 Put also $\,L(i,p,k)=\sL_i/p^k\sL_i\,$, where $\,k\,$ is a positive
 integer, $\,p>3\,$ is a prime number, $\,i\in\set{1,2,4,6}\,$;
 $\,\La^2(p,k)=\La^2/p^k\La^2\,$ and $\,\rS^2(p,k)=\rS^2/p^k\rS^2\,$,
 where $\,p\ge3\,$ is again a prime number.

\begin{theorem}
\begin{enumerate}
 \item
 The string and band cubic modules defined above, the modules
 $\,L(i,p,k)\,$ and the modules $\,\La^2(p,k),\,\rS^2(p,k)\,$ are just
 all finitely generated indecomposable 2-divisible cubic modules.  
 \item
  The only isomorphisms between these modules are the following:
 \begin{itemize}
  \item
  $\,M^\rD\iso M^{\rD^*}$, where $\,\rD^*\,$ denotes the diagram
 symmetric to $\,\rD\,$, the latter being of type \emph{(ii)} or
 \emph{(iii)}; 
  \item
  $\,M^\rB\iso M^{\rB^{(s)}}\,$, where
 $\,\rB^{(s)}=(\rD^{(s)},f(t))\,$ and $\,\rD^{(s)}\,$ is the
 \emph{$\,s$-th shift} of the diagram $\,\rD\,$, i.e., the diagram
$$
  \begin{array}{*{11}c}
 j_{2s+1}& & & &j_{2s+2}j_{2s+3}&  & & &j_{2s} \\
 &\kk{k_{2s+1}}& &\kk{k_{2s+2}}& &\kk{\cdots}&
 &\kk{k_{2s}}& \\ 
 & & i_{2s+1}i_{2s+2}& & & &i_{2s-1}i_{2s}& &
 \end{array}
$$
  \item
  $\,M^\rB\iso M^{{\rB^*}^{(s)}}$, where
  $\,(\rD,f(t))^*=(\rD^*,f^*(t))\,$ 
 and $\,f^*(t)=\la_1^{-1}t^df(1/t)\,$;
 \end{itemize}
\end{enumerate}
\end{theorem}
\begin{proof}
[\rm Cf.] \cite[Theorem 5.1]{qu}
\end{proof}
\begin{corol}
 \begin{enumerate}
 \item A 2-divisible cubic module $\,M\,$ is of finite
  projective dimension \iff it contains no string summands of types
  {\rm (i)} and {\rm (iii)}. In this case either $\,M\,$ is projective
  (hence, a direct sum of modules isomorphic to $\,\sB_i\,$) or
  $\,\prd M=1\,$. Otherwise, $\,\prd M=\8\,$.
  \item Every 2-divisible cubic module $\,M\,$ has a periodic projective
 resolution 
$$
\begin{CD}
 \dots\to P_{n+1}\gai nP_n@>\ga_{n-1}>>\dots\gai1P_1\gai0P_0\to M\to0
\end{CD}
$$
 of period 6, namely, such that $\,\ga_{n+6}=\ga_n\,$ for every
 $\,n\ge2\,$.  
 \item
  If $\,M\,$ is an indecomposable non-projective 2-divisible cubic
 module, $\,\nT M\,$ is its torsion part, then 
$$
 M/\nT M\iso \begin{cases}
	\sL_{i_{2n}} &\text{\rm if }\,M \text{ \rm is a string module
 	of type (i)}\\
	\sL_{i_1}\+\sL_{i_{2n}} &\text{\rm if }\,M \text{ \rm is a
	string module of type (ii)}\\
	0 &\text{\rm otherwise}
	     \end{cases}
$$
\rm (Remind that here $\,i_{2n}\,$ denotes the unique index such that
 $\,i_{2n-1}\sim i_{2n}\,$ and $\,i_1\,$ denotes the unique index such
 that $\,i_2\sim i_1\,$.)  
\end{enumerate}
\end{corol}
\begin{proof}
[\rm Cf] \cite[Corollaries 5.3, 5.5]{qu}.
\end{proof}

\section{Weakly alternative functors}
\label{sec4}

 Consider now the case of cubic $\,\mZ$-modules $\,F\,$ such that
 $\,F(\mZ)=0\,$. We call them \emph{weakly alternative}. They are
 actually modules over the category $\,\bA^a=\bA/(\id_1)\,$. The
 functors $\,\La^3\,$ and $\,\La^2\*\Id\,$ are of this sort. Here
 again, we are able to obtain a complete description.

 Denote by $\,\bC_0\,$ the subring in $\,\mZ\xx\Mat(2,\mZ)\,$ consisting
 of the pairs $\,(a,B)\,$ with $\,b_{12}\equiv0,\ a\equiv
 b_{11}\pmod3\,$. 

\begin{prop}
 The category $\,\bA^a\,$ is Morita equivalent to the semi-direct
 product $\,\bC=(\mZ\xx\bC_0)\ltimes T\,$ where $\,T\,$ is the
 elementary $\,2$-group with three generators $\,\xi,\eta,\th\,$ such
 that the following relations hold:
\begin{align*}
& \th=\xi\eta,\ \eta\xi=0,\\
& \eta e=\eta,\ e\xi=\xi,\ \text{ \rm where }\ e=(1,0,0),\\
& \xi\left(0,a,B\right)=a\xi\quad \text{\rm and}\quad
 \left(0,a,B\right)\eta=a\eta 
\end{align*}
 for any pair $\,(a,B)\in\bC_0\,$.
\end{prop}
\begin{proof}
 The category $\,\bA^a\,$ has two objects: $\,2,3\,$, and is generated
 by the morphisms $\,h_i:2\to3,\,p_i:3\to2\,$ (i=1,2) subject to the
 following relations, cf. diagram \eqref{eq2}:
\begin{align*}
& h_ip_ih_i=2h_i\,,\quad p_ih_ip_i=2p_i\,,\quad h_ip_j=0\,\textrm{ if
}\, i\ne j\,,\\
& h_1+h_2=h_1p_1h_2p_2h_1+h_2p_2h_1p_1h_2 \,,\\
& p_1+p_2=p_1h_2p_2h_1p_1+p_2h_1p_1h_2p_2 \,,
\end{align*}
 that imply:
\begin{align*}
 h_ip_i&=h_ip_ih_jp_jh_ip_i\,,\\
 2p_i&=2p_ih_jp_jh_ip_i\,,\\
 2h_i&=2h_ip_ih_jp_jh_i
\end{align*}
 (always $\,i,j\in\set{1,2},\,i\ne j\,$). Then $\,e_i=p_ih_jp_jh_i\,$
 are orthogonal idempotents in $\,\bA^a(2,2)\,$ which are conjugate as
 $\,e_1p_1h_2=p_1h_2e_2\,$ and $\,e_2p_2h_1=p_2h_1e_1\,$. Moreover,
 one can check that $\,e_3p_1h_2=p_1h_2e_3=e_3p_2h_1=p_2h_1e_3\,$
 where $\,e_3=1-e_1-e_2\,$. Denote this product by  $\,\th\,$. Then
 $\,2\th=0\,$ and $\,p_ih_i=\th+2e_i\,$, so $\,\bA^a(2,2)\,$ is Morita
 equivalent to $\,\mZ\xx\mZ[\,\th\,]\,$.

 Just in the same way, one can verify that
 $\,\bA^a(3,3)=\gnr{f_1,f_2,\al_1,\al_2,\be}\,$ where
 $\,f_1=h_1p_1h_2p_2,\ f_2=1-f_1,\ \al_i=h_ip_i,\
 \be=\al_2\al_1\,$. Moreover, $\,f_i\al_i=\al_if_j=\al_i\,$,
 $\,\al_1\al_2=3f_1\,$ and $\,\be^2=3\be\,$. Therefore,
 $\,\bA^a(3,3)\,$ is isomorphic to the subring $\,\bC_0\,$.

 Let $\,u=p_1h_2p_2,\,v=h_1p_1h_2p_2h_1\,$. Then $\,e_1u=uf_1,\,
 f_1v=ve_1,\, uv=e_1,\,vu=f_1\,$, so $\,e_1\,$ and $\,f_1\,$ are
 conjugate in $\,\bA^a\,$. Put $\,\xi=p_1-p_1h_2p_2h_1p_1,\,
 \eta=h_1-h_1p_1h_2p_2h_1\,$. One can also check that
\begin{align*}
& 2\xi=0,\ 2\eta=0,\ \eta\xi=0,\ \xi\eta=\th,\
 \xi f_3=\xi,\ f_3\eta=\eta,\\& e_3\bA^a(3,2)=\gnr{\xi},\
 \bA^a(2,3)e_3=\gnr{\eta}\,. 
\end{align*}
 Hence, putting $\,E^a_3=\bA^a(2,\_)e_3\,$
 we get that $\,\bA^a\,$ is Morita equivalent to the endomorphism ring
 $\,\End(E_3^a\+\bA^a(3,\_))\,$ which is isomorphic to $\,\bC\,$.
\end{proof}

 Describe now (finitely generated) $\,\bC$-modules or, the same,
 weakly alternative cubic functors. If $\,M\,$ is such a module, put
 $\,\oM=M/TM\,$. It is a module over $\,\mZ\xx\bC\,$, so
 $\,\oM=M_1\+M_2\,$, where $\,M_1\,$ is an abelian group with
 $\,\bC_0M_1=0\,$ and $\,M_2\,$ is a $\,\bC_0$-module. Again, the same
 observations as in \cite[Theorem 2.1]{qu} yield the following result.
\begin{prop}
\label{locwa}
 \begin{enumerate}
  \item 
  Weakly alternative cubic modules $\,M,N\,$ are isomorphic \iff
  $\,M_{(p)}\iso N_{(p)}\,$ for every prime $\,p\,$.   
  \item 
  Given a set $\,\set{N_p}\,$ where $\,N_p\,$ is a weakly
  alternative cubic $\,\mZ_{(p)}$-module, there is a weakly
  alternative cubic module  $\,M\,$ such that $\,M_{(p)}\iso N_p\,$
  for all $\,p\,$ \iff almost all of them
  are torsion free (maybe, zero) and $\,\mQ\*N_p\iso\mQ\*N_q\,$ for
  all $\,p,q\,$. In this case $\,\mQ\*M=\mQ\*N_p\,$ and $\,\nT
  M=\bop_p\nT N_p\,$.
 \end{enumerate}
\end{prop}

 Let $\,\oC=\bC/T=\mZ\xx\bC_0\,$, $\,\sC_1=\sL_1=\oC e\,$,
 $\,\sC_2=\oC(0,1,e_{11})\,$ and $\,\sC_3=\sL_4=\oC(0,0,e_{22})\,$.
 They are all indecomposable projective $\,\oC$-modules. Put also
 $\,\sL_2=\oC(0,1,0)\,$ and $\,\sL_3=\oC(0,0,e_{11})\,$. They are all
 non-projective indecomposable torsion free $\,\oC$-modules. If
 $\,p>2\,$, the localization $\,\bC_{(p)}\,$ is torsion free. If,
 moreover, $\,p>3\,$, this localization is hereditary, hence, any
 $\,\bC_{(p)}$-module 
 splits into a direct sum of a torsion free and a torsion one, the
 latter being a direct sum of modules isomorphic to
 $\,L(i,k,p)=\sL_i/p^k\,$ for $\,i\in\set{1,2,3},\,k\in\mN\,$. 

 The description of $\,\bC_{(3)}$-modules is similar to that of
 $\,\bB$-modules in the preceding section. One only has to consider
 the set $\,\set{1,2,3,4}\,$ instead of $\,\set{1,2,3,4,5,6}\,$ with
 the relations $\,-\,$ and $\,\sim\,$ defined as follows: $\,3-4\,$
 and $\,2\sim3\,$. 

 At last, $\,\bC_{(2)}=\bC'\xx\bC''\,$ where
 $\,\bC'=(\mZ_{(2)}\xx\mZ_{(2)})\ltimes T\,$ and
 $\,\bC''\iso\Mat(2,\mZ_{(2)})\,$. A $\,\bC''$-module is a direct sum
 of several copies of $\,(\sL_3)_{(2)}\,$ and of modules isomorphic to
 $\,L(3,k,2)=\sL_3/2^k\,$. A $\,\bC'$-module $\,W\,$ is given by a
 diagram of $\,\mZ_{(2)}$-modules of the form:
$$
\begin{CD}
	W_1 \mat{ @>\xi>> \\ @<<\eta< } W_2
\end{CD}
$$
 such that $\,2\xi=0,\,2\eta=0,\,\eta\xi=0\,$. Split both $\,W_1\,$
 and $\,W_2\,$ into a direct sum of cyclic modules $\,C_i=\mZ/2^i\,$
 and $\,C_\8=\mZ_{(2)}\,$. Then one can check that such
 a diagram is a direct sum of diagrams $\,W(\om)\,$ and
 $\,W(\om,\pi)\,$. Here $\,\om\,$ is a (finite) word of the form 
$$
  \dots\xi^{i_m}\eta_{j_m}\xi^{i_{m+1}}\eta_{m+1}\xi\dots \quad
  (i_l,j_l\in\mN\cup\set{\8})
$$
 containing no subwords of the form $\,_\8\xi,\ ^\8\eta,\
 \eta_1\xi\,$. In $\,W(\om,\pi)\,$, $\,\om\,$ must be of the form
$$
  _j\xi^{i_1}\eta_{j_1}\xi^{i_2}\eta_{j_2}\xi\dots^{i_m}\eta_j
$$
 with the same restrictions as above
 and $\,\pi(t)\ne t^n\,$ is a primary polynomial over $\,\mZ/2\,$.
 Namely, if $\,W=W(\om)\,$, then $\,W_1=\bop_mC_{i_m}\,$,
 $\,W_2=\bop_mC_{j_m}\,$, while $\,\xi(C_{i_m})\sbe C_{j_{m-1}}\,$,
 $\,\eta(C_{j_m})\sbe C_{i_m}\,$ and the induced mappings are the
 unique non-zero ones of period 2 (we denote them by $\,\ga\,$). If
 $\,W=W(\om,\pi)\,$ and $\,\deg\pi=n\,$, then $\,W_1=\bop_mnC_{i_m}\,$,
 $\,W_2=nC_j\+(\bop_mnC_{j_m})\,$; $\,\xi(nC_{i_m})\sbe nC_{j_{m-1}}\,$,
 $\,\eta(nC_{j_m})\sbe nC_{i_m}\,$, $\,\xi(nC_{i_1})\sbe nC_j\,$,
 $\,\eta(nC_j)\sbe nC_{i_m}\,$ and the induced mappings are given by
 the matrices $\,\ga I\,$, except the last one, given by
 $\,\ga\Phi\,$ where $\,\Phi\,$ is the Frobenius cell corresponding
 to the polynomial $\,\pi\,$. 

 Note that the torsion free part of $\,W(\om,\pi)\,$ is zero, that of
 $\,W(\om)\,$ consists of at most one cyclic summand, and that of an
 indecomposable $\,\bC_{(3)}$-module either is trivial, or consists of
 one or of two cyclic summand (for string modules $\,M^\rD\,$ of type,
 respectively, (i)  or (ii), cf. page  \pageref{str1}). So, in
 accordance with Proposition~\ref{locwa}, the indecomposable
 $\,\bC$-modules $\,M\,$ that are not torsion have the following local
 components $\,M_{(p)}\,$ (we only describe $\,M_{(2)}\,$ and
 $\,M_{(3)}\,$ as all other ones are torsion free, hence, uniquely
 defined):
\begin{enumerate}
  \item
  $\,M_{(2)}=W(\om)\,$ where $\,\om\,$ contains $\,\xi^\8\,$,
  $\,M_{(3)}=M^\rD\,$  
  where $\,\rD\,$ is a string of type (i) or (ii) with
  $\,i_{2n-1}=2\,$ or $\,i_2=2\,$ (if both  
  $\,i_2=i_{2n-1}=2\,$ in a string of type (ii), there are two
  possibilities for such $\,M\,$). 
  \item
  $\,M_{(2)}=W(\om)\+W(\om')\,$ where both $\,\om\,$ and $\,\om'\,$
  contain $\,\xi^\8\,$,  
  $\,M_{(3)}=M^\rD\,$ where $\,\rD\,$ is a string of type (ii) with
  $\,i_2=i_{2n-1}=2\,$. 
  \item
  $\,M_{(3)}\,$ is a string module of type (i) or (ii), $\,M_{(2)}\,$
  is torsion free (thus, uniquely defined).
\end{enumerate}
 (Note that the case when $\,M_{(3)}\,$ is torsion free is a part of
 case (1) above.) 
 
\section{Cubic vector spaces}
\label{sec5}

 Now we consider the \emph{cubic vector spaces}, i.e., cubic functors
 $\,F:\fab\to\mod{}\fK\,$ where $\,\fK\,$ is a field. If
 $\,\chr\fK\ne2\,$ they are a partial case of the functors considered
 in Section~\ref{sec3}; hence, we always suppose that
 $\,\chr\fK=2\,$. In this case a cubic functor $\,F\,$ is given by a
 diagram of $\,\fK$-vector spaces of the same shape \eqref{eq2}
 with the relations:
\begin{equation}\label{eq1a}
\begin{split}
 &  h_1p_2=h_2p_1=0\,,\quad h_1h=h_2h\,,\quad pp_1=pp_2\,,\\
 &  h_ip_ih_i=0\,,\quad p_ih_ip_i=0\quad (i=1,2)\,,\\
 &  hph=0\,,\quad php=0\,,\\
 &  \ol hp+h_1+h_2=h_1p_1h_2p_2h_1+h_2p_2h_1p_1h_2\,,\\
 &  h\ol p+p_1+p_2=p_1h_2p_2h_1p_1+p_2h_1p_1h_2p_2
\end{split}
\end{equation}
 (just as before, $\,\ol h=h_1h=h_2h\,$ and $\,\ol
 p=pp_1=pp_2\,$). Denote by $\,\sC\,$ the $\,\fK$-linear category with
 objects $\,1,2,3\,$ and generating morphisms $\,h:1\to2,\ p:2\to 1,\
 h_i:2\to 3,\ p_i:3\to 1\,$ ($\,i=1,2\,$) satisfying the relations
 \eqref{eq1a}. So, a cubic vector space is the same as a
 $\,\sC$-module (i.e., a linear functor $\,\sC\to\mod{}\fK\,$).

 The last two equations imply that
$$
 \ol h\ol p+h_ip_i=h_ip_ih_jp_jh_ip_j\quad (i,j=1,2;\ i\ne j),
$$
 whence $\,p_i\ol h\ol p=\ol h\ol ph_i=0\,$ for $\,i=1,2\,$. Hence,
 the elements $\,e_i=h_ip_ih_jp_j\,$ as well
 as the elements $\,f_i=p_ih_jp_jh_i\,$ ($\,i,j=1,2;\ i\ne j\,$) are
 orthogonal idempotents, respectively, in $\,\sC(3,3)\,$ and in
 $\,\sC(2,2)\,$. Thus, in $\,\add\sC\,$, $\,3\iso3_0\+3_1\+3_2\,$ and
 $\,2\iso2_0\+2_1\+2_2\,$, so that the identity morphisms of $\,3_i\,$
 are identified with $\,e_i\,$ (with $\,e_0=1-e_1-e_2\,$) and those of
 $\,2_i\,$ are identified with $\,f_i\,$ (with $\,f_0=1-f_1-f_2\,$).
 In what follows, we write $\,\sC(x,y)\,$ for the set of morphisms
 $\,x\to y\,$ in $\,\add\sC\,$. An easy calculation shows that the
 four objects $\,2_i,\,3_i\,$ ($\,i=1,2\,$) are isomorphic in
 $\,\add\sC\,$. For instance, as
 $\,p_1h_2p_2=p_1h_2p_2h_1p_1h_2p_2\,$, this element lies in
 $\,\sC(3_1,2_1)\,$; the element $\,h_1p_1h_2p_2h_1\,$ lies in
 $\,\sC(2_1,3_1)\,$ and their products are just $\,f_1\,$ and
 $\,e_1\,$, whence $\,2_1\iso3_1\,$, etc. So, we only have to take into 
 account one of these objects, say $\,3_1\,$. One can also easily
 check that $\,\sC(3_1,3_1)=\fK e_1\,$, while
 $\,\sC(x,3_1)=\sC(3_1,x)=0\,$ if $\,x\in \set{1,2_0,3_0}\,$. Thus,
 the category $\,\sC\,$ is Morita equivalent to the direct product of
 the \emph{trivial category} with one object $\,3_1\,$ and the full
 subcategory $\,\sC^*\,$ of $\,\add\sC\,$ with the objects
 $\,1,2_0,3_0\,$. One easily check that the cubic module $\,T\,$
 corresponding to the (unique) indecomposable representation of the
 trivial part is the following one:
\begin{equation}\label{triv}
\begin{split}
&  T_1=0\,,\quad T_2=\gnr{u_1,u_2}\,,\quad T_3=\gnr{v_1,v_2}\,;\\
&  h_1(u_1)=v_1\,,\ h_1(u_2)=0\,;\quad h_2(u_1)=0\,,\ h_2(u_2)=v_2\,;\\ 
&  p_1(v_1)=0\,,\ p_1(v_2)=u_1\,;\quad p_2(v_1)=u_2\,,\ p_2(v_2)=0\,.
\end{split}
\end{equation}
 (This cubic module corresponds to the functor $\,\fK\*\La^2\*\Id\,$.) 
 So, from now on, we only consider the representations of $\,\sC^*\,$
 and for every morphism $\,\al\,$ from $\,\sC\,$ we denote by the same
 letter $\,\al\,$ its restriction onto $\,\sC^*\,$. As such
 restrictions of $\,e_i\,$ and $\,f_i\,$ are zero for $\,i=1,2\,$, one
 obtains the relations:
$$
  h\ol p+p_1+p_2=0\,,\qquad \ol hp+h_1+h_2=0\,,
$$
 so we may exclude $\,p_2,h_2\,$ from the generating set. Therefore,
 $\,\sC^*$-modules are just diagrams of vector spaces
\begin{equation}\label{eq2b}
\begin{CD}
F_1 \mat{@>h>>\\@<p<<} F_2 \mat{@>{h_1}>>\\@<{p_1}<<} F_3 
\end{CD}
\end{equation}
 with the relations
\begin{equation}\label{eq1b}
 hph=php=h_1p_1h_1=p_1h_1p_1=0\,,\quad h_1p_1=h_1hpp_1\,.
\end{equation}
 (One easily checks that they imply all relations \eqref{eq1a} if we
 put $\,h_2=h_1+h_1hp\,$ and $\,p_2=p_1+hpp_1\,$.)

 Consider the subdiagram
$$
\begin{CD}
F_1 \mat{@>h>>\\@<p<<} F_2 \,.
\end{CD}
$$
 As $\,hph=php=0\,$, it decomposes into a direct sum of
 the following shape:
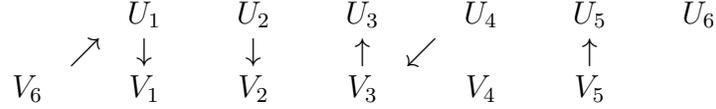
\begin{equation*}\def\hh{\hspace*{.3cm}}
\begin{array}{*{13}c}
  &\hh&U_1&\hh&U_2&\hh&U_3&\hh&U_4&\hh&U_5&\hh&U_6\\
     &\nearrow&\downarrow&&\downarrow&&\uparrow&\swarrow&&&\uparrow&&\\
     V_6& &V_1& &V_2& &V_3& &V_4& &V_5& &
\end{array}
\end{equation*}
 ($\,U_i\,$ are the direct summands of $\,F_1\,$, $\,V_i\,$ those of
 $\,F_2\,$, the arrows show the action of $\,h\,$ and $\,p\,$ when it
 is non-zero, the corresponding maps being isomorphisms.) It means
 that $\,h\,$ and $\,p\,$ are given by the following matrices:
\begin{equation}\label{F12}
h=\mtr{I&0&0&0&0&0\\ 0&I&0&0&0&0\\ 0&0&0&I&0&0&\\
       0&0&0&0&0&0\\ 0&0&0&0&0&0\\ 0&0&0&0&0&0 }\,,\quad
p=\mtr{0&0&0&0&0&I\\ 0&0&0&0&0&0\\ 0&0&I&0&0&0\\
       0&0&0&0&0&0\\ 0&0&0&0&I&0\\ 0&0&0&0&0&0}\,.
\end{equation}
 ($\,I\,$ denotes the identity matrices.)

 With respect to the decomposition of $\,F_2\,$, $\,h_1\,$ and
 $\,p_1\,$ are given by the matrices
 $\,H=(H_1,H_2,H_3,H_4,H_5,H_6)\,$ and
 $\,P=(P_1,P_2,P_3,P_4,P_5,P_6)^\top\,$  where $\,H_i:V_i\to F_3,\
 P_i:F_3\to V_i\,$. For these matrices the following conditions
 hold: 
\begin{itemize}
 \item
 the number of rows of $\,H\,$ equals the number of columns of
 $\,P\,$;
 \item
  the number of columns of $\,H_i\,$ is the same as the number of rows
 of $\,P_i\,$ for every $\,i\,$; 
 \item
  the number of columns of $\,H_1\,$ is the same as the number of rows
 of $\,H_6\,$.
\end{itemize}

 When one applies the isomorphisms of the spaces $\,F_i\,$ which do
 not destroy the form \eqref{F12} of $\,h\,$ and 
 $\,p\,$, they result in elementary transformations of the columns of
 the matrices $\,H\,$ and $\,P\,$ such that:
\begin{itemize}
 \item
  the transformations of $\,P\,$ are \emph{contragredient} to those of
  $\,H\,$ (e.g., when we add the $\,k$-th column of $\,H\,$ to the
  $\,l$-th one, we have to subtract the $\,l$-th column of $\,P\,$
  from the $\,k$-th one, etc.);
 \item
  the transformations inside $\,H_1\,$ are the same as those inside
  $\,H_6\,$;
 \item
  one can only add columns of $\,H_i\,$ to those of $\,H_j\,$ if
  $\,i\le j\,$ and $\,(i,j)\ne(3,4)\,$.
\end{itemize}
  Using such transformations, one can reduce the matrix $\,H\,$ to the
  following form:
$$
\left(
\begin{array}{cccc|cc|ccc|ccc|cc|cccc}
 0&0&I&0&0&0&0&0&0&0&0&0&0&0&0&0&0&0\\
 0&0&0&I&0&0&0&0&0&0&0&0&0&0&0&0&0&0\\
 0&0&0&0&0&I&0&0&0&0&0&0&0&0&0&0&0&0\\
 0&0&0&0&0&0&0&I&0&0&I&0&0&0&0&0&0&0\\
 0&0&0&0&0&0&0&0&I&0&0&0&0&0&0&0&0&0\\
 0&0&0&0&0&0&0&0&0&0&0&I&0&0&0&0&0&0\\
 0&0&0&0&0&0&0&0&0&0&0&0&0&I&0&0&0&0\\
 0&0&0&0&0&0&0&0&0&0&0&0&0&0&0&I&0&0\\
 0&0&0&0&0&0&0&0&0&0&0&0&0&0&0&0&0&I\\
 0&0&0&0&0&0&0&0&0&0&0&0&0&0&0&0&0&0
\end{array}
\right)
$$
 (the vertical lines denote the borders of the matrices $\,H_i\,$).
 Certainly, we have to subdivide in the same manner the rows and
 columns of the matrix $\,P\,$. Denote the corresponding blocks by
 $\,P_{ij}\,$ ($\,1\le i\le17,\,1\le j\le10\,$); the horizontal
 (vertical) stripes of $\,P\,$ will be denoted by $\,R_i\,$
 (respectively, $\,S_j\,$). Then the relations
 \eqref{eq1b} for $\,h_1\,$ and $\,p_1\,$ are equivalent to the
 following conditions for these blocks: 
\begin{equation*}
\begin{split}
& P_{ij}=0\ \text{ if }\ i\in\set{4,6,9,12,14,16,18},\\
& P_{ij}=0\ \text{ if }\ i\in\set{3,17}\ \text{ and }\ j\ne\,10\,,\\ 
& P_{8j}=P_{11,j}\ \text{ for all }\ j\,,\\
& P_{3,10}=P_{17,10}\ \text{ and }\ P_{i1}P_{3,10}=0\ \text{ for all
}\ i\,. 
\end{split}
\end{equation*}
 Thus, we only have to consider the matrix $\,\oP\,$ obtained from
 $\,P\,$ by crossing out all zero horizontal stripes  as well as the
 stripes $\,R_3\,$ and $\,R_8\,$.
 A straightforward calculation shows that the automorphisms of the
 spaces $\,F_3\,$ and $\,V_i\,$ which do not destroy the shape of
 $\,H\,$ give rise to the elementary transformations of the matrix
 $\,\oP\,$ satisfying the following conditions:
\begin{itemize}
 \item
  the transformations inside $\,R_1\,$ are the same as those inside
  $\,R_{15}\,$; 
 \item
  the transformations inside $\,S_2\,$ are the same as those inside
  $\,S_9\,$;
 \item
  the transformations inside $\,R_2\,$ are contragredient to those
  inside $\,S_8\,$;
 \item
  the transformations inside $\,R_{11}\,$ are contragredient to those
  inside $\,S_4\,$; 
 \item
  the transformations inside $\,R_{17}\,$ are contragredient to those
  inside $\,S_1\,$;
 \item
  one can add the columns of $\,S_i\,$ to those of $\,S_j\,$ \iff
  $\,i\le j\,$ and $\,(i,j)\ne(5,6)\,$;
 \item
  one can add the rows of $\,R_i\,$ to those of $\,R_j\,$ \iff
  $\,i\ge j\,$ and $\,(i,j)\ne(10,7)\,$.
\end{itemize}
 Obviously, one may suppose that the matrix $\,P_{17,10}\,$ is of the
 form $\,\mtr{0&I\\0&0}\,$. Then the whole stripe $\,S_1\,$ has the
 form $\,\mtr{0&S_1'}\,$, the number of columns in $\,S_1'\,$ being the
 same as that of zero rows in $\,P_{17,10}\,$. Moreover, using the
 elementary transformations described above, one may suppose that the
 remaining part of the stripe $\,S_{10}\,$ is of the form
 $\,\mtr{S'_{10}&0}\,$, the number of columns in $\,S'_{10}\,$ being
 the same as in the zero part of $\,P_{17,10}\,$. In what follows, we
 omit the dash and write $\,S_1\,$, $\,S_{10}\,$, $\,P_{i1}\,$ and
 $\,P_{i,10}\,$ for the remaining (non-zero) parts of the corresponding
 matrices. Of course, we should no more consider the stripe
 $\,R_{17}\,$. 

 Now one immediately sees that we are again in the situation
 considered in \cite{bo}. Namely, we have got the following
 semi-chains: 
\begin{equation*}
\begin{split}
 \dE&=\setsuch{R_i}{R_1>R_2>R_5>R_7>R_{11}>R_{13}>R_{15},\
 R_5>R_{10}>R_{11}}\,,\\ 
 \dF&=\setsuch{S_j}{S_1<S_2<S_3<S_4<S_5<S_7<S_8<S_9<S_{10},\,
 S_4<S_6<S_7}
\end{split}
\end{equation*}
 with the involution $\,\si\,$ such that 
$$
 \si(R_1)=R_{15}\,,\ \si(R_2)=S_8\,,\ \si(R_{11})=S_4\,,\ \si(S_2)=S_9 
$$
 and $\,\si(x)=x\,$ if $\,x\notin\set{R_1,R_2,R_{11},R_{15},S_2,S_4,S_8,S_9}\,$. 
% with the equivalence relation $\,\sim\,$ described as follows:
% $$
% R_1\sim R_{15}\,,\ R_2\sim S_8\,,\ R_{11}\sim S_4\,,\ S_2\sim S_9\,.
% $$
 Hence, one can deduce from \cite{bo} a list of canonical forms of the
 matrix $\,P\,$ and, therefore, of cubic functors
 $\,\fab\to\mod{}\fK\,$. It can again be arranged in the form of
 ``strings and bands,'' though a trifle more sophisticated than in the
 preceding sections.

\begin{defin}
 Put $\,\dX=(\dE\cup\dF)\=\set{R_{10},S_6}\,$. Write $\,x-y\,$ if
 $\,x\in\dE,\,y\in\dF\,$ or vice versa; $\,x\sim y\,$ if
 $\,\si(x)=y\ne x\,$ or $\,x=y\in\set{R_7,S_5}\,$. Call
 $\,R_7,\,S_5\,$ \emph{special elements}. An
 $\,\dX$-\emph{word} (or simply \emph{word}) is a sequence
 $\,w=x_1r_2x_2r_3\dots r_nx_n\,$, where $\,x_k\in\dX\,$,
 $\,r_k\in\set{\sim,-}\,$, satisfying the following conditions: 
 \begin{itemize}
  \item
  $\,x_{k-1}r_kx_k\,$ (in the above defined sense) for all
  $\,k=2,\dots,n\,$; 
  \item
  $\,r_k\ne r_{k-1}\,$ for all $\,k=3,\dots,n\,$;
 \end{itemize}
 Such a word is called \emph{full} if the following conditions hold:
 \begin{itemize}
  \item
  either $\,r_2={\sim}\,$ or $\,x_1\not\sim y\,$ for every $\,y\ne
 x_1\,$;
  \item
  either $\,r_n={\sim}\,$ or $\,x_n\not\sim y\,$ for every $\,y\ne x_n\,$.
 \end{itemize}
 A full word is called \emph{special} if $\,x_1\,$ or $\,x_n\,$, but
 not both of them, is a special element, \emph{bispecial} if both
 $\,x_1\,$ and $\,x_n\,$ are special, and \emph{ordinary} if neither
 $\,x_1\,$ nor $\,x_n\,$ is special. A word $\,w\,$ is called
 \emph{non-symmetric} if $\,w\ne w^*\,$, where $\,w^*\,$ is the
 \emph{inverse word}: $\,w^*=x_nr_nx_{n-1}\dots x_2r_2x_1\,$.

 An $\,\dX$-word is called \emph{cyclic} if $\,r_2=r_n=-\,$ and
 $\,x_n\sim x_1\,$. A cyclic word is called \emph{aperiodic} if it is
 not of the form $\,v\sim v\sim \dots\sim v\,$ for a shorter word
 $\,v\,$. The $\,s$-th \emph{shift} of a cyclic word $\,w\,$ is the
 word $\,w^{(s)}=x_{2s+1}r_{2s+2}x_{2s+2}\dots x_n\sim x_1r_2\dots
 r_{2s}x_{2s}\,$. A cyclic word $\,w\,$ is said to be
 \emph{shift-symmetric} if $\,w^{(s)}\,$ is symmetric for some
 $\,s\,$. Note that the length $\,n\,$ of a cyclic word is always
 even.
\end{defin}

\begin{defin}
 A \emph{string datum} $\,\rD\,$ is:
\begin{itemize}
  \item
  either an ordinary non-symmetric $\,\dX$-word $\,w\,$ (\emph{ordinary
  string datum});
  \item
  or a pair $\,(w,\de)\,$ consisting of a special word $\,w\,$ and
  $\,\de\in\set{0,1}\,$ (\emph{special string datum});
  \item
  or a quadruple $\,(w,\de_1,\de_2,m)\,$ consisting of a bispecial
  non-sym\-metric word $\,w\,$, $\,\de_1,\de_2\in\set{0,1}\,$ and
  $\,m\in\mN\,$ (\emph{bispecial string datum}).
\end{itemize}
 Put $\,\rD^*=w^*\,$ in the first case, $\,\rD^*=(w^*,\de)\,$ in the
 second and $\,\rD^*=(w^*,\de_2,\de_1,m)\,$ in the third one.

 A \emph{band datum $\,\rB\,$} is a pair $\,(w,\pi(t))\,$ consisting
 of an aperiodic cyclic word $\,w\,$ and of a \emph{primary
 polynomial} $\,\pi(t)\in\fK[t]\,$ (i.e., a power of an irreducible
 one) such that $\,\pi(t)\ne t^d\,$ if $\,w\,$ is not shift-symmetric
 and $\,\pi(t)\ne(t-1)^d\,$ if $\,w\,$ is shift-symmetric. Put
 $\,\rB^{(s)}=(w^{(s)},f(t))\,$ and
 $\,\rB^*=(w^*,\la^{-1}t^d\pi(1/t))\,$, where $\,w\,$ is not
 shift-symmetric, $\,d=\deg\pi\,$ and $\,\la=\pi(0)\,$. 
\end{defin}

 Now the results of \cite{bo} immediately imply the following
\begin{theorem}
 \begin{enumerate}
  \item
  Every string data $\,\rD\,$ defines an indecomposable cubic vector
  space $\,V^\rD\,$ called \emph{string cubic space}.
  \item
  Every band data $\,\rB\,$ defines an indecomposable cubic vector
  space $\,V^\rB\,$ called \emph{band cubic space}.
  \item
  Every indecomposable cubic vector space, except
  $\,\fK\*\La^2\*\Id\,$, is isomorphic either to a string or to a band
  one. 
  \item
  The only isomorphisms between string and band spaces are the
  following: 
  \begin{itemize}
   \item
   $\,V^\rD\iso V^{\rD^*}\,$ where $\,\rD\,$ is a string data;
   \item
   $\,V^\rB\iso V^{\rB^{(s)}}\,$ and $\,V^\rB\iso
  V^{{\rB^*}^{(s)}}\,$ where $\,\rB\,$ is a band data and
  $\,s\in\mN\,$. 
  \end{itemize}
 \end{enumerate}
\end{theorem}
 Moreover, one can deduce from \cite{bo} and the preceding
 considerations an explicit construction of string and band
 spaces, though it is rather cumbersome and we will not include it
 here. 

\section{Torsion free cubic modules}
\label{sec6}

 In this section, we consider \emph{torsion free} cubic modules, i.e.,
 cubic functors $\,\fab\to\fab\,$. Again, we first study them
 \emph{locally}, i.e., describe cubic functors
 $\,\fab\to\fab_{(p)}\,$, the latter being the category of torsion
 free finitely generated $\,\mZ_{(p)}$-modules. Denote also by
 $\,\mZ_p\,$ the ring of $\,p$-adic integers. The latter has the
 advantage of being complete, that guarantees lifting idempotent
 endomorphisms modulo $\,p\,$. Note that the calculations of
 Section~\ref{sec3} easily imply that
 $\,\bA\*\mQ\iso\mQ^3\xx\Mat(2,\mQ)\xx\Mat(4,\mQ)^2\,$. In particular,
 it is a semi-simple \emph{split} $\,\mQ$-algebra, i.e., $\,\End
 V\iso\mQ\,$ for every simple module $\,V\,$. Thus, it follows from
 the standard results of the theory of lattices over orders
 (cf.~\cite[Chapter 4, \S\,1]{rog}) that every torsion free
 $\,\bA_p$-module is actually a completion of a torsion free
 $\,\bA_{(p)}$-module. Hence, lifting idempotents is possible for
 $\,\bA_{(p)}$-modules too. Together with the calculations of
 Section~\ref{sec5} it implies the following result.
\begin{prop}
\label{loc2}
 The ring $\,\bA_{(2)}\,$ is isomorphic to the subring of
 $\,\mZ_{(2)}^3\xx\Mat(2,\mZ_{(2)})\xx\Mat(4,\mZ_{(2)})^2\,$
 consisting of  all sextuples $\,(a_1,a_2,a_3,B,C,D)\,$ satisfying the
 following congruences \rm modulo $2\,$:
\begin{align*}
 & a_1\equiv b_{11}\equiv c_{11}\,,\ a_2\equiv b_{22}\equiv
 c_{22}\equiv c_{33}\,,\ a_3\equiv c_{44}\,,\\
 & b_{12}\equiv0\ \text{ and }\ c_{ij}\equiv0\ \text{ if }\
 i<j\,. \tag{*}
\end{align*}
\end{prop}

 One can easily check that this ring is an example of the so called
 \emph{Backstr\"om order} \cite{bac}, i.e., its Jacobson
 radical coincides with that of an hereditary order
 $\,\bH\,$. Thus, torsion free $\,\bA_{(2)}$-modules are in a natural
 \oc with the representations of a quiver $\,\sQ\,$ over the field
 $\,\mZ/2\,$. Namely, the vertices of $\,\sQ\,$ are just simple
 $\,\bA_{(2)}$-modules $\,A_1,\dots,A_r\,$ and simple $\,\bH$-modules
 $\,H_1,\dots,H_s\,$, all arrows are from some $\,A_i\,$ to some
 $\,H_j\,$ and the number of such arrows equals the multiplicity of
 $\,A_i\,$ in $\,H_j\,$. In our example $\,\bH\,$ consists of the
 sextuples satisfying the congruences (*) only, $\,r=4,\,s=10\,$ and
 the quiver $\,\sQ\,$ consists of $\,4\,$ connected components:
$$
 \begin{array}{cp{1cm}c}
 \begin{CD}
 H_1 @<<< A_1 @>>> H_4 \\
     &&  @VVV   & \\
  &&    H_6 & 
 \end{CD}& &
 \begin{CD}
  && H_5 & \\
  && @AAA & \\
  H_2 @<<< A_2 @>>> H_7 \\
  & & @VVV & \\
  && H_8 &
 \end{CD}\\&&\\
\begin{CD}  H_3 @<<< A_3 @>>> H_9 \end{CD} &&
\begin{CD} A_4 @>>> H_{10} \end{CD}
 \end{array}
$$
 (The numeration is the natural one with respect to the description of
 $\,\bA_{(2)}\,$ and $\,\bH\,$ above.) This quiver is tame and the
 list of its representations is well known (cf., e.g., \cite{dlr}).
 Hence, we can derive the description of torsion free
 $\,\bA_{(2)}$-modules.
 
 The description of indecomposable torsion free $\,\bA_{(3)}$-modules
 is given in Section~\ref{sec3}. For any other prime $\,p\,$, the
 localization $\,\bA_{(p)}\,$ is just
 $\,\mZ_{(p)}^3\xx\Mat(2,\mZ_{(p)})\xx\Mat(4,\mZ_{(p)})^2\,$, hence,
 there are exactly $6$ indecomposable (and irreducible) torsion free
 modules. Therefore, the standard ``local--global'' procedure
 \cite{rog,ad} implies the following result on torsion free cubic
 modules. 

\begin{prop}
 Torsion free cubic modules are in \oc with the pairs
 $\,(M_2,M_3)\,$, where $\,M_p\,$ is a torsion free
 $\,\bA_{(p)}\,$-module ($\,p=2,3\,$) and $\,\mQ\*M_2\iso\mQ\*M_3\,$. 
\end{prop}
\begin{proof}
 It follows from \cite[Chapter 4]{rog} that such a pair $\,(M_2,M_3)\,$ always
 defines a cubic module $\,M\,$ \emph{up to genus}. (Remind that two modules
 $\,M,\,N\,$ \emph{belong to the same genus} if $\,M_p\iso N_p\,$ for
 all prime $\,p\,$.) Note that
 $\,\Ga=\mZ^3\xx\Mat(2,\mZ)\xx\Mat(4,\mZ)^2\,$ is a
 maximal order containing $\,\bA\,$ and two torsion free
 $\,\Ga$-modules belonging to the same genus are isomorphic. Applying
 the results of \cite{ad}, we see that the isomorphism classes of
 modules belonging to the same genus as $\,M\,$ are in \oc with the
 double cosets
$$
  \Aut(\Ga M)\backslash \Aut(\Ga M_2)\xx\Aut(\Ga M_3)/ \Aut(M_2)\xx
  \Aut(M_3) \,,
$$
 where $\,\Ga M\,$ denotes the $\,\Ga$-submodule in $\,\mQ\*M\,$
 generated by $\,M\,$. As $\,\bA\sp6\Ga\,$, we can replace these
 cosets by 
$$
 \ol{\Aut}(\Ga M)\backslash \Aut(\Ga M/6\Ga M)/\Aut(M/6\Ga M)\,,
$$
 where $\,\ol{\Aut}(\Ga M)\,$ denotes the image of $\,\Aut(\Ga M)\,$
 in $\,\Aut(\Ga M/6\Ga M)\,$. But $\,\End(\Ga M)\,$ is
 just a direct product of matrix algebras $\,\Mat(n_i,\mZ)\,$, and any
 matrix invertible modulo $6$ is the image of an invertible integer
 matrix. Hence, $\,\ol{\Aut}(\Ga M))= \Aut(\Ga M/6\Ga M)\,$, so every
 genus only contains one module up to isomorphism.
\end{proof}

 Using the arguments analogous to those of \cite[Theorem 2.1]{qu}, one
 gets the following corollary  extending Theorem~\ref{t21} 
 to \emph{all} cubic functors.
\begin{corol}
\begin{enumerate}
 \item 
  Cubic modules $\,M,N\,$ are isomorphic \iff
 $\,M_{(p)}\iso N_{(p)}\,$ for every  prime $\,p\,$.   
\item 
 Given a set $\,\set{N_p}\,$ where $\,N_p\,$ is a
 cubic $\,\mZ_{(p)}$-module, there is a cubic module
  $\,M\,$ such that $\,M_{(p)}\iso N_p\,$ for all $\,p\,$ \iff almost
  all of them are torsion free (maybe,
  zero) and $\,\mQ\*N_p\iso\mQ\*N_q\,$ for all $\,p,q\,$. In this case
  $\,\mQ\*M=\mQ\*N_p\,$ and $\,\nT M=\bop_p\nT N_p\,$.
\end{enumerate}
\end{corol}
 It seems plausible that the analogous result is no more true for the
 polynomial functors of degree $\,4\,$, but at the moment we do not
 have a corresponding counterexample.

\section{One conjecture}
\label{sec7}

 The descriptions of quadratic modules and of cubic \emph2-divisible modules
 as well as some other evidence give rise to the following conjecture
 concerning polynomial modules of higher degrees.

 Put $\,\mZ_{<p}=\mZ[\,1/(p-1)!\,]\,$ and denote by $\,\bA^{(p)}\,$ the
 subring of the direct product
 $\,\mZ_{<p}\xx\mZ_{<p}\xx\Mat(2,\mZ_{<p})^{p-1}\,$ consisting of all
 $\,(p+1)$-tuples 
$$
    \left(a,b,C^1,\dots,C^{p-1}\right)
$$
 such that $\,\md a{c^1_{11}}p,\ \md b{c^{p-1}_{22}}p\,$ and
 $\,\md{c^i_{22}}{c^{i+1}_{11}}p\,$ for every $\,i=1,\dots,p-2\,$.
\begin{conj}
 The category $\,\kM(p)=\mod p{\mZ_{<p}}\,$ is equivalent to the
 category of modules over $\,\bA^{(p)}\xx\mZ_{<p}^r$ for appropriate
 $\,r\,$ (depending on $\,p\,$).
\end{conj}

 As the ring $\,\bA^{(p)}\,$ fits the conditions of \cite{pure}, it
 would give a complete description of such polynomial
 functors in strings and bands terms quite similar to that of
 Section~\ref{sec3}. Moreover, easy calculations show that this
 conjecture would imply the following properties of such polynomial
 functors analogous to those of quadratic and \emph2-divisible cubic
 functors: 
\begin{itemize}
 \item
  Two modules from $\,\kM(p)\,$ are isomorphic \iff all their
  localizations are isomorphic. 
 \item
  Given a set $\,\set{N_q}\,$ of polynomial functors ($N_q\in\mod
  p{\mZ_{(q)}}$, $\,q\,$ runs through all primes $\,\ge p\,$), there is
  a functor $\,M\in\kM(p)\,$ such that $\,M_{(q)}\iso N_q\,$ for all
  $\,q\,$ \iff almost all $\,N_q\,$ are torsion free and
  $\,\mQ\*N_q\iso \mQ\*N_{q'}\,$ for all $\,q,q'\,$. In this case
  $\,\nT M=\bop_q\nT N_q\,$ and $\,\mQ\*M\iso\mQ\*N_q\,$.
 \item
  Any functor from $\,\kM(p)\,$ is either projective, or of projective
  dimension 1, or of infinite projective dimension.
 \item
  Any functor from $\,\kM(p)\,$ has a periodic projective resolution
  of period $\,2p\,$.
 \item
  If a functor $\,F\in\kM(p)\,$ is indecomposable and non-projective,
  its torsion free part is a direct sum of at most two irreducible
  torsion free modules.
\end{itemize}

\section*{Acknowledgments}

 This research was mainly performed when the author was at the
 Max-Planck-Institut f\"ur Mathematik. I am extremely grateful to the
 Institute for the excellent opportunity and working conditions. I am also
 grateful to H.-J.~Baues who enthusiastically supported this research
 and to W.~Dreckman for useful discussions.

\end{document}